\def\zbar{{\bar{z}}}

\def\wbar{{\bar{w}}}

\def\ID{{\mathbb D}}

\def\S2{{\bf S}(2)}

\def\IR{{\mathbb R}}

\def\IC{\mathbb C}

\def\1ton{1,2,\ldots,n}

\def\div{{\rm div}}

\long\def\symbolfootnote[#1]#2{\begingroup\def\thefootnote{\fnsymbol{footnote}}
\footnote[#1]{#2}\endgroup}

\documentclass[a4paper]{article}

\usepackage{amsthm, amssymb}
\usepackage{amsfonts}
\usepackage{mathrsfs}
\usepackage{epsfig,multicol}
\usepackage{makeidx}
\usepackage{graphicx}
\usepackage{textcase}
\usepackage{setspace} 

\makeindex

\newtheorem{theorem}{Theorem}[section]

\newtheorem{corollary}[theorem]{Corollary}

\newtheorem{conjecture}[theorem]{Conjecture}

 \title{Homeomorphic Solutions to 
 Reduced \\ Beltrami Equations }
 \author{Kari Astala \and Jarmo J\"a\"askel\"ainen \thanks{Authors were supported by the Academy of Finland, project
no.\ 1118634, the Finnish Centre of Excellence in
Analysis and Dynamics Research, and  project MRTN-CT-2006-035651, Acronym
CODY, of the European Commission.}}
 
\begin{document}
 \date{}

 \maketitle


\frenchspacing
 
\begin{abstract}\noindent We  study differential expressions related to linear families of quasiconformal mappings and give a simple and direct proof to a result due to Alessandrini and Nesi \cite{ANess}.
\end{abstract}

\vspace{0.5cm}

{\footnotesize
\noindent \textit{2000 AMS Mathematics Classification Numbers.} Primary 30C62. Secondary 35J45.

}

\section{Introduction}

The reduced Beltrami differential equation
\begin{equation}\label{yksi}
  \frac{\partial f}{\partial \zbar} = \lambda(z) \; \Im m  \left( \frac{\partial f}{\partial z} \right)\!\!,
   \qquad |\lambda(z)| \leqslant k < 1, 
\end{equation}
arises naturally in many different contexts in the theory of quasiconformal mappings; for instance in the Stoilow factorization and the $G$-closure problems for the general Beltrami equation
\begin{equation}\label{kaksi}
  \frac{\partial f}{\partial \zbar} = \mu(z) \frac{\partial f}{\partial z} + \nu(z)\overline{\,\frac{\partial f}{\partial z}\,}, \qquad |\mu(z)| + |\nu(z)|   \leqslant k < 1.
\end{equation}
The solutions to (\ref{yksi}) have a number of special properties, such as $f(z) \equiv z$ being always a solution. This is the unique normalized solution fixing  $0$ and $1$. In fact, 
if  $g \in W^{1,2}_{loc}(\IC)$ is a homeomorphism satisfying the equation  (\ref{yksi}) and  $g$ has two fixed points, then $g(z) \equiv z$.  We refer the reader for these and other properties of the reduced equation to the monograph \cite{AIM}. An early application of the reduced equation  for uniqueness properties of (\ref{kaksi}) for self-mappings of the unit disk can be found in  \cite{Boj}.

Studies of the reduced Beltrami equation (\ref{yksi}) indicate that  for its solutions, the  null Lagrangian
$$
{\cal J}(z,f) =  \Im m  \left( \frac{\partial f}{\partial z} \right)
$$
has many  properties analogous to the familiar Jacobian determinant of a Sobolev mapping. This suggest the following conjecture, cf. \cite[p. 222]{AIM}.

\begin{conjecture} \label{conje2}
Suppose $f: \Omega \to \IC$ is a (quasiregular) solution to  the reduced Beltrami equation {\rm (\ref{yksi})}.
Then  either $ \Im m \big(f_z\big)$  is a  constant, or else
$$
 \Im m \left( \frac{\partial f}{\partial z} \right)  \not= 0 \qquad \mbox{almost everywhere in $\Omega$.}
$$
\end{conjecture}
\medskip
\noindent Through the Stoilow factorization type theorems, and \cite[Theorem 6.1.1]{AIM} in particular, the conjecture  can equivalently be formulated in terms of  solutions to the general  Beltrami system (\ref{kaksi}), see Section 2.

In fact, for homeomorphic solutions the conjecture is closely related to the notion \cite{AIM, BDIS} of linear  families  of quasiconformal mappings. 
Given a domain $\Omega \subset \IC$ and an $\IR$-linear subspace ${\cal F} \subset W^{1,2}_{loc}(\Omega)$, we say  that   ${\cal F} $ is  a {\it linear family of quasiconformal mappings},  if  there is  $1 \leqslant K< \infty$ such that for every $g \in  {\cal F}$, either $g \equiv 0$ or else $g$ is a $K$-quasiconformal mapping in $\Omega$. 
It quickly follows \cite{BDIS} that ${\rm dim}\, {\cal F} \leqslant 2$. If we have the equality, then
$$ {\cal F} =\{ a\, \Phi + b\, \Psi: \; a, b \in \IR  \}
$$
for some quasiconformal mappings $ \Phi: \Omega \to \IC$ and $ \Psi: \Omega \to \IC$. In this case we say that the family ${\cal F} $ is {\it generated} by the mappings $\Phi$ and $\Psi$.
In particular,  if $\Phi$ and $\Psi$ generate a linear family of quasiconformal mappings, then by definition each $F_{a, b} = a\, \Phi + b\, \Psi$ is injective in $\Omega$, whenever $a^2 + b^2 \neq 0$.

In general,   quasiconformality is not preserved under linear combinations. However, if we have  mappings that happen to be  solutions to the same Beltrami equation (\ref{kaksi}), then their linear  combinations are at least quasiregular. Conversely,  \cite{BDIS} associates to a linear (two-dimensional) family ${\cal F} $ of quasiconformal mappings a Beltrami equation of the type (\ref{kaksi}) satisfied by every $g \in {\cal F}$. The next theorem, answering in positive \cite[Conjecture 1]{BDIS},
 implies that the associated equation is unique. 

\smallskip

\begin{theorem}\label{conje} Let ${\cal F}$ be a linear family of quasiconformal mappings in a domain  $\Omega \subset \IC$. If $\Phi, \Psi \in {\cal F}$,  then either 
$$
{\cal J}(\Phi, \Psi) :=  \; \Im m  \left( \frac{\partial \Phi}{\partial z}  \overline{\, \frac{\partial \Psi}{\partial z}\,} \,\right) \neq 0 \qquad \mbox{almost everywhere in $\Omega$}
$$
or else
$$ a\,  \Phi(z) + b \, \Psi (z) \equiv 0 \qquad \mbox{for some constants }  a, b \in \IR, \; a^2 + b^2 \neq 0,$$
in which case  ${\cal J}(\Phi, \Psi) \equiv 0$.
\end{theorem}
\smallskip

\noindent It is possible, in fact,  to obtain 
this theorem by combining results and methods from 
\cite{ANess} and \cite{BDIS}. However, the purpose of this paper is to give a simple and direct  proof to this beautiful result.
 Our methods in proving the theorem are similar to those of Alessandrini and Nesi in \cite{ANess}, but we will simplify their  approach.  For principal solutions the  result was also announced in \cite{Boj2}, but unfortunately \cite[Proposition 2]{Boj2} is not  valid, with  counterexamples  easy to find.\medskip

Not every pair of homeomorphic solutions generate a linear family of quasiconformal mappings; simple examples can be found already among  the solutions to the Cauchy-Riemann system $f_\zbar =0$. For instance, $\Phi(z) = z$ and $\Psi(z) = z^2$ are  both conformal in $\Omega = \{z: \Re e(z) >0 \}$, yet (some of) their linear combinations are non-injective in $\Omega$, and thus the mappings do not generate a linear family of (quasi)conformal mappings. 

However, global homeomorphic solutions to (\ref{kaksi}) are determined by their values at two distinct points, and it follows from this fact that  in the domain $\Omega= \IC$, linear combinations of  homeomorphic solutions are either constants or homeomorphisms, see \cite[Section 6.2]{AIM}. Hence  Theorem \ref{conje} applies.
\vskip9pt

\begin{corollary}\label{viisi}
Suppose $\Phi, \Psi \in W^{1,2}_{loc}(\IC)$ are homeomorphic solutions to {\rm (\ref{kaksi})}. Then, unless the mappings are affine combinations of each other, 
$$\; \Im m  \left( \frac{\partial \Phi}{\partial z}  \overline{\, \frac{\partial \Psi}{\partial z} \,} \,\right) \neq 0 \qquad \mbox{almost everywhere in $\IC$.}
$$
\end{corollary}
\medskip

\noindent As an immediate consequence, Conjecture \ref{conje2} follows for global homeomorphic solutions $f:\IC \to \IC$ to the reduced equation (\ref{yksi}).
\medskip

There are further situations where the injectivity of a linear family of  solutions to (\ref{kaksi}) can be guaranteed. For example,  if  $f$ satisfies (\ref{kaksi}) in a bounded convex domain $\Omega$ with
$$ \Re e \big(f(z)\big) =  \Re e \big({\cal A}(z)\big) \quad \mbox{on}\;  \partial \Omega, \qquad 
{\cal A}:\IC \to \IC \;\; \mbox{a linear isomorphism},
$$
 then $f$ is injective: With the Stoilow factorization we can write $f = h \circ F^{-1}$, where $h$ is holomorphic in the unit disk $\ID$ and $F:\ID \to \Omega$ is a quasiconformal homeomorphism. Since ${\cal A} \circ F $ maps to a convex domain, by  the classical Rad\'o-Kneser-Choquet theorem the Poisson extension $U$ of its boundary values is one-to-one. But $ \Re e \big(h\big) =  \Re e \big(U\big)$ on  $\partial \ID$, hence in $\ID$, and by the theorem of Clunie and Sheil-Small \cite[Theorem 5.3]{Clunie}, or \cite[p. 38]{Duren},  the injectivity of  $h$ and $f$ follows. 
For alternative proofs of injectivity, using  properties of the Beltrami equation, see  \cite{ANess}, \cite{BDIS} or \cite{Nesileonetti}. We now obtain the following result   of  Alessandrini and Nesi in   \cite{ANess}.

\begin{corollary}\label{seiska} Suppose $\Omega \subset \IC$ is a bounded convex domain, and let 
$\Phi, \Psi \in W^{1,2}(\Omega)$ be solutions to {\rm (\ref{kaksi})} in $\Omega$, such that
$$ \Re e \big( \Phi(z) - z \big) \in  W^{1,2}_{0}(\Omega), \qquad \Re e \big(\Psi(z) + i z \big) \in  W^{1,2}_{0}(\Omega).
$$
Then 
$$\; \Im m  \left( \frac{\partial \Phi}{\partial z}  \overline{\, \frac{\partial \Psi}{\partial z}\,} \,\right) > 0 \qquad \mbox{almost everywhere in $\Omega$.}
$$
\end{corollary}
\medskip

\noindent According to \cite[Lemma 7.1]{BDIS} we have $\Im m  \big( \Phi_ z  \overline{\Psi_ z}\, \big) \geqslant 0$ almost everywhere. Hence the  corollary  follows from Theorem \ref{conje}.
\medskip

Finally we note that the quantity ${\cal J}(\Phi, \Psi)$  arises naturally  in the study  \cite{BDIS, GIKMS}  of the  $G$-closure problems of Beltrami operators, and this connection was also the motivation in the work of Alessandrini and Nesi. In fact, combining Theorem \ref{conje}  with  results and ideas developed in  \cite{BDIS} and \cite{GIKMS}, we see that the family ${\cal F}_K(\IC)$, $1 \leqslant K < \infty$,  of  Beltrami differential operators
\[
\frac{\partial }{ \partial \overline{z}}- \mu(z) \frac{\partial }
{ \partial z}-\nu(z) \overline{\frac{\partial }{ \partial z}}
\]
with
\[
|\mu(z)|+|\nu(z)|   \leqslant  \frac{K-1 }{ K+1} = k < 1, \qquad z \in \IC,\]
is $G$-compact. For the details, we refer the reader to 
\cite[Chapter 16]{AIM}.

\section{Proof of Theorem 1.2}

We start by reducing  Theorem \ref{conje}  to the reduced Beltrami equation, using basic facts from \cite{AIM} and \cite{BDIS}. For this we may assume that  $\Phi, \Psi \in W^{1,2}_{loc}(\Omega)$ generate the linear family ${\cal F}$ of $K$-quasiconformal mappings. 
According to \cite[Section 5.3]{BDIS} or \cite[Remark 16.6.7]{AIM}, there are Beltrami coefficients 
$\mu$ and $\nu$ such that
 every element $g \in {\cal F}$ satisfies the equation
\begin{equation}\label{yleinen}
\frac{\partial g}{\partial \zbar} = \mu(z) \frac{\partial g}{\partial z} + \nu(z)\overline{\,\frac{\partial g}{\partial z}\,} \qquad \mbox{almost everywhere in $\Omega$},
\end{equation}
where
$$ |\mu(z)| + |\nu(z)| \leqslant \frac{K-1}{K+1} = k < 1.
$$

Next, following \cite[Lemma 7.1]{BDIS}, we apply the fact that $\Phi$ and  
$\Psi$ generate a linear family of injections. Since for every $a, b\in \IR$ the mappings $ a\,  \Phi(z) + b \, \Psi (z)$ are injections, we have 
$$\Lambda(z,w):= \Im m \left( \frac{\Phi(z) - \Phi(w)}{\Psi(z)-\Psi(w)}\right) \neq 0, \qquad z, w \in  \Omega, \; z \neq w.
$$
As the complement of the diagonal $\{(z,z): z \in \Omega \}$ is connected  in $\Omega \times \Omega$, the continuous function $\Lambda(z,w)$ does not change sign. We may assume that $\Lambda(z,w) < 0$ whenever $z \neq w$; otherwise replace $\Psi$ by $-\Psi$.
From this fact and Taylor's first-order development, we obtain at points $z$ of differentiability, thus almost everywhere, that
$$
\Im m  \left(  \Phi_z(z) \overline{\,  \Psi_z(z)\,} \,\right) \leqslant 0.
$$
The explicit  details can be found  in \cite[Lemma 7.1]{BDIS} and in \cite[p. 203]{AIM}.

We now arrive at a reduced equation. Namely,
$$ \Psi (z) = (f \circ \Phi)(z), \qquad z \in \Omega,
$$
for some quasiconformal homeomorphism $f:\Phi(\Omega) \to \Psi(\Omega)$.
The general Stoilow factorization  theorem \cite[Theorem 6.1.1]{AIM} states that, since $\Phi$ and $\Psi$ satisfy the same equation (\ref{yleinen}), the mapping $f$ is a solution to the reduced equation,  
$$ \frac{\partial f}{\partial \wbar} = \lambda(w) \; \Im m  \left( \frac{\partial f}{\partial w} \right)\!\!, \qquad w \in \Phi(\Omega).
$$
Here $\lambda(w) = -2 i \nu(z)/(1+|\nu(z)|^2 - |\mu(z)|^2)$ and $w = \Phi(z)$. Furthermore, by  the ellipticity bounds in (\ref{yleinen}), $|\lambda(w)| \leqslant k' = 2k/(1+k^2) < 1$. 

With the chain rule one calculates  that
\begin{equation}
\label{kasi}
J(z, \Phi)\, \Im m \big( f_w \circ \Phi \big) = (-1 +|\mu|^2 -|\nu|^2 )\,  \Im m \big( \Phi_z \overline{\Psi_z}\, \big) \geqslant 0
\end{equation}
almost everywhere. In particular, as quasiconformal mappings preserve Lebesgue zero sets, $ \Im m \big( f_w  \big) \geqslant 0$ almost everywhere in $\Omega' = \Phi(\Omega)$. Moreover,
 Theorem \ref{conje} is equivalent to showing that $ \Im m \big( f_w  \big)$ can vanish only in a set of measure zero.
\vskip8pt

With this reduction, we are now left to study the homeomorphic  solution $f \in W^{1,2}_{loc}(\Omega')$ to the reduced equation (\ref{yksi}). Let us  write $f(z) = u(z) + iv(z)$, where $u$ and $v$ are real valued. Similarly write $\lambda(z) = \alpha(z) + i\beta(z)$. 
 
 Taking  the imaginary part of (\ref{yksi}) shows us that
$u_y + v_x = \beta(v_x-u_y)$,  i.e.
\[  u_y = \frac{\beta-1}{\beta+1} \, v_x. \]
Thus 
\begin{equation}\label{viimet457}
 2 \, \Im m \left( \frac{\partial f}{\partial z} \right) = v_x - u_y =  \frac{2}{\beta+1} \, v_x = \frac{2}{\beta-1} \,u_y.
\end{equation}
Since $|\beta(z)| \leqslant |\lambda(z)| \leqslant k < 1$, the coefficients $2/(\beta(z)\pm 1)$ in  (\ref{viimet457}) are uniformly bounded below. Hence to prove Theorem \ref{conje} it suffices to show that $u_{y} \neq 0$  almost everywhere. 
\medskip

The trick of the proof is that, for the reduced equation (\ref{yksi}), the derivative  $u_y$ is a solution to  the adjoint equation determined by  a non-divergence type operator.  To state this more precisely, consider an operator
$$
L = \sum_{i,j = 1}^{2} \sigma_{ij}(z)\,\frac{\partial^{2}}{\partial x_{i}x_{j}},
$$
where $\sigma_{ij} = \sigma_{ji}$ are measurable and the matrix
$$
\sigma(z) = 
\left[ \begin{array}{cc}
\sigma_{11}(z) & \sigma_{12}(z)  \\
\sigma_{12}(z) & \sigma_{22}(z)  
\end{array} \right]
$$
is uniformly elliptic, 
$$
\frac{1}{K}|\xi|^{2} \leqslant \langle\sigma(z)\xi, \xi\rangle 
= \sigma_{11}(z)\xi_{1}^{2} + 2\sigma_{12}(z)\xi_{1}\xi_{2} + \sigma_{22}(z)\xi_{2}^{2} \leqslant K|\xi|^{2}
$$
for all $ \xi\in\IC$ and $z \in \Omega'$. Here $K$ is the ellipticity constant.
Then we say that the  function $w\in L^{2}_{{loc}}(\Omega')$ is a weak solution to the adjoint 
equation  
\begin{equation}\label{viimet342}
L^{*}w = 0
\end{equation}
if
$$
\int wL\varphi = 0, \qquad \mbox{for every $\varphi\in C^{\infty}_{0}(\Omega')$.} 
$$

To identify $u_y$ as a solution to an equation of the type (\ref{viimet342}), we recall that  the components of solutions $f=u + iv\;$ to general Beltrami equations satisfy a divergence type second-order equation, see \cite[Section 16.1.5]{AIM}. In case of (\ref{yksi}), it turns out that the component $u$ satisfies the equation
\begin{equation}\label{div}
\div\,  A\nabla u = 0, \qquad
A(z) := 
\left[ \begin{array}{cc}
1 & a_{12}(z)  \\
0 & a_{22}(z)  
\end{array} \right]\!\!,
\end{equation}
where the matrix elements are
\begin{equation}\label{viimet333}
a_{12} = \frac{2\, \Re e (\lambda)}{ 1- \Im m(\lambda)}\, = \frac{2\, \alpha}{1-\beta}, \quad 
a_{22} = \frac{1 + \Im m(\lambda)}{ 1- \Im m(\lambda)}  = \frac{1+\beta}{1-\beta} > 0.
\end{equation}
Precisely, (\ref{div}) means that  for every $\varphi\in C^{\infty}_{0}(\Omega')$ we have 
\begin{equation}\label{ch6id}
0 = \int \nabla\varphi \cdot A\nabla u = \int \varphi_{x} (u_{x} + a_{12}u_{y}) + \varphi_{y}a_{22}u_{y}.
\end{equation}
But since derivatives of smooth test functions are again test functions, we can replace   $\varphi$ by  $\varphi_{y}\in C^{\infty}_{0}(\Omega')$. In this case  the identity (\ref{ch6id}) takes  the form
\begin{eqnarray*}
0 &=&  \int \varphi_{yx}u_{x} + a_{12}\varphi_{yx}u_{y} + a_{22}\varphi_{yy}u_{y} \\
& = & \int \varphi_{xx} u_{y} + a_{12}\varphi_{xy}u_{y} + a_{22}\varphi_{yy}u_{y} \\
&= & \int (\varphi_{xx}  + a_{12}\varphi_{xy} + a_{22}\varphi_{yy})u_{y}.
\end{eqnarray*}
Thus $u_{y}$ is indeed a distributional solution to the adjoint 
equation $L^{*}u_{y} = 0$, where
\begin{equation} \label{viimet345}
L = \frac{\partial^{2}}{\partial x^{2}} + a_{12} \, \frac{\partial^{2}}{\partial x\partial y} 
+ a_{22}\, \frac{\partial^{2}}{\partial y^{2}}\,
\end{equation}
and $a_{12}, a_{22}$ are given by (\ref{viimet333}).
Note that the original matrix $A(z)$ is not symmetric. However,  the operator $L$ in  (\ref{viimet345}) can  be represented by the symmetric matrix 
$$
\sigma(z) := 
\left[ \begin{array}{cc}
1 & a_{12}(z)/2  \\
a_{12}(z)/2 & a_{22}(z)  
\end{array} \right]\,
$$
and as $|\lambda(z)| \leqslant k' < 1$, from (\ref{viimet333}) we see that $\sigma$ is uniformly elliptic.
\medskip

Next,  we use (\ref{kasi}) and (\ref{viimet457}) to prove that the derivative $u_y \geqslant 0$ almost everywhere.  In fact, it is precisely here  we use the assumption that  ${\cal F}$ consists only of homeomorphisms.
\smallskip

Thus we may assume that $u_y$ is a non-negative solution to the adjoint equation, $L^*u_y = 0$,  where $L$ is defined in (\ref{viimet345}). In this case we may simply apply a  result of Fabes and Stroock \cite{fab} that the non-negative solutions to the adjoint equation satisfy a uniform reverse H\"older estimate.  
For the reader's convenience we recall here the explicit formulation of their theorem.

\begin{theorem} {\rm \cite[Theorem 2.1]{fab}\label{rh2}} \label{rh345} Consider the operator
$$
L = \sum_{i,j = 1}^{2} \sigma_{ij}(z)\,\frac{\partial^{2}}{\partial x_{i}x_{j}},
$$
where $\sigma_{ij} = \sigma_{ji}$ are measurable and $L$ is uniformly elliptic with constant $K$. 

Then there  exists a constant $C_0$, depending only on the ellipticity constant $K$,
such that for all $w \geqslant 0$ satisfying $L^{*}w = 0 $ in a domain $\Omega \subset \IC$ we have 
 \begin{equation}\label{viimet708}
\biggl[\frac{1}{r^{2}}\int_{\ID(z_0,r)}w(z)^{2} \,dz \biggr]^{1/2} 
\leqslant   \frac{C_0}{r^{2}}\int_{\ID(z_0,r)}w(z) \, dz
\end{equation}
in every disk $\,\ID(z_0,r)$ such that  $\ID(z_0,2r)\subset\Omega$.
 \end{theorem}
\medskip

Applying the Fabes-Stroock theorem to $u_y = w$, it follows from  (\ref{viimet708}) that either $w \equiv 0$ or $w > 0$ almost everywhere. Namely,
 if  $E:= \{z \in \Omega' : w(z) = 0\}$ has positive measure and  $z_0 \in E$ is a point of density, we can find disks $\ID_r = \ID(z_0,r)$ with 
\[  |\ID_r \setminus E| < \varepsilon r^2, \qquad \mbox{ where } C_0 \sqrt{\varepsilon} \, < \, 1,
\]
and $C_0$ is the constant of the reverse 
H\"older inequality (\ref{viimet708}). Thus
\begin{equation} \label{viimet009}
 \int_{\ID_r}  w = \int_{\ID_r\setminus E}  
w \leqslant\left( \int_{\ID_r} w^2\right)^{1/2} |\ID_r \setminus E|^{1/2}  
\leqslant C_0 \sqrt{\varepsilon}\,  \int_{\ID_r}  w.
\end{equation}
This is possible only if  $w$ vanishes identically in $\ID(z_0,r)$, that is $\ID(z_0,r) \subset E$. 
But then  we can replace the disk by a slightly larger one, argue as in (\ref{viimet009}) and by iterating the argument prove that $E = \Omega'$.  Now $u_y \equiv 0$, and  (\ref{yksi}) with (\ref{viimet457}) shows that $f$ is holomorphic with the real valued derivative, hence  affine. Therefore the mappings $\Phi$ and $\Psi = f \circ \Phi$  would not generate a linear family of homeomorphisms. 
Thus $u_y\neq 0$ almost everywhere, and we have  completed the proof of Theorem \ref{conje}. \hfill $\Box$

\bigskip
 \bigskip
 
 kari.astala@helsinki.fi 
 
 jarmo.jaaskelainen@helsinki.fi
 
 \end{document}